\documentclass{amsart}

\usepackage{amsmath,amstext,amsthm,amscd,amsopn,verbatim,amssymb}
\usepackage{hyperref}
\usepackage{mathbbol}
\usepackage{graphicx}
\usepackage{color}
\usepackage{pstricks, pstricks-add, pst-node, pst-coil}
\usepackage{multirow}
\usepackage{array}

\newtheorem{thm}{Theorem}
\newtheorem{prop}[thm]{Proposition}
\newtheorem{lemma}[thm]{Lemma}
\theoremstyle{definition}

\newtheorem{cor}[thm]{Corollary}

\newtheorem{question}[thm]{Question}

\theoremstyle{remark}

\newcommand{\alg}[1]{\mathfrak{{#1}}}
\newcommand{\tr}{\text{tr}}

\newcommand{\co}[2]{\left[{#1},{#2}\right]} 

\newcommand{\eref}[1]{(\ref{#1})} 

\newcommand{\ad}{{\text{ad}}}

\newcommand{\p}{\partial}

\newcommand{\C}{{\mathbb{C}}}

\DeclareMathOperator{\li}{Li}

\title{On the (ir)rationality of Kontsevich weights}
\author{Giovanni Felder}
\email{felder@math.ethz.ch}

\author{Thomas Willwacher}
\email{thomas.willwacher@math.ethz.ch}

\address{Department of Mathematics, ETH Zurich}
\thanks{This work was partially supported by the Swiss National Science Foundation (grant 200020-105450)}
\date{}
\begin{document}

\begin{abstract}
We compute the weight of a Kontsevich graph in deformation quantization. Up to rationals, the result is $\zeta(3)^2/\pi^6$.
\end{abstract}
\maketitle

\section{Introduction and the main result}
There are three open questions on the rationality of Kontsevich integrals in deformation quantization:
\begin{question}[strong version]
\label{conj:strong}
Is the weight of any Kontsevich graph, obtained using Kontsevich's harmonic propagator, rational?
\end{question}
\begin{question}[strong Lie version]
\label{conj:stronglie}
Is the weight of any Lie graph, obtained using Kontsevich's harmonic propagator, rational?
\end{question}
\begin{question}[weak version]
Is the universal star product, obtained using Kontsevich's harmonic propagator, rational?
\end{question}

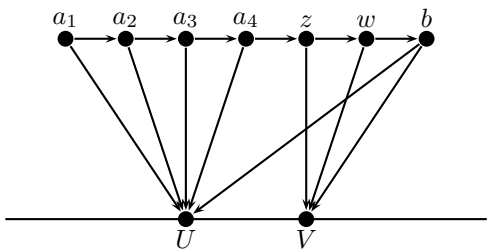
\begin{figure}[t]
\psset{unit=.8cm}
\psset{arrowscale=1}
\begin{pspicture}(0,-1)(8,4)  
  \cnode*(1,3){3pt}{a1}
  \cnode*(2,3){3pt}{a2}
  \cnode*(3,3){3pt}{a3}
  \cnode*(4,3){3pt}{a4}

  \cnode*(5,3){3pt}{z}  
  \cnode*(6,3){3pt}{w}
  \cnode*(7,3){3pt}{b}

  \cnode*(3,0){3pt}{l}  
  \cnode*(5,0){3pt}{r}
  
  \psline{-}(0,0)(8,0)
  \psset{arrows=->}
  \ncline{a1}{a2}
  \ncline{a2}{a3}
  \ncline{a3}{a4}
  \ncline{a4}{z}
  \ncline{z}{w}
  \ncline{w}{b}
  
  \ncline{a1}{l}
  \ncline{a2}{l}
  \ncline{a3}{l}
  \ncline{a4}{l}
  \ncline{b}{l}
  \ncline{z}{r}
  \ncline{w}{r}
  \ncline{b}{r}

  \uput{1ex}[-90](l){$U$}
  \uput{1ex}[-90](r){$V$}  

  \uput{1ex}[90](a1){$a_1$}
  \uput{1ex}[90](a2){$a_2$}  
  \uput{1ex}[90](a3){$a_3$}
  \uput{1ex}[90](a4){$a_4$}  
  \uput{1ex}[90](z){$z$}
  \uput{1ex}[90](w){$w$}  
  \uput{1ex}[90](b){$b$}
\end{pspicture}
\caption{\label{fig:theLiegraph} This Lie graph is the subject of Theorem \ref{thm:main}.}
\end{figure}

Note that the ``weak Lie version'' that one is tempted to formulate is actually a well-known Theorem, see Appendix \ref{sec:kontscbh}. The result of this paper is the following:
\begin{thm}
\label{thm:main}
The Kontsevich weight of the Lie graph in Figure \ref{fig:theLiegraph} is, up to rationals, $\frac{\zeta(3)^2}{\pi^6}$.
\end{thm}

It is a famous open problem in number theory to determine whether $\frac{\zeta(3)}{\pi^3}$ is algebraic or even rational. The Theorem shows that from a positive answer to Question \ref{conj:strong} or \ref{conj:stronglie} algebraicity would follow. 

%
%

\section{Preliminaries and definitions}
\subsection{Graphs and weights}
A \emph{Kontsevich graph} \cite{kontsevich} is a directed graph with two types of vertices, type I and type II. The type I vertices are usually denoted $1,\dots,m$ ($m\geq 0$) and the type II vertices $\bar{1},\dots,\bar{n}$ ($n\geq 0$). The directed edges are required to start at type I vertices. In fact, we will only consider the case $m=7$, $n=2$. A \emph{Lie graph} is a Kontsevich graph with exactly two type II vertices, and with at most one edge ending and exactly two edges starting at every type I vertex. An example is shown in Figure \ref{fig:theLiegraph}.

To each Kontsevich graph, one can associate a \emph{weight}, i.e., a real, conjecturally rational number. The weight is given by an integral of the form
\[
\int_{C} \bigwedge_{(u,v)}\nu(u,v).
\]
Here the domain of the integration is essentially the space of all embeddings of the vertex set of the graph $\Gamma$ into the closure of the Poincar\'e (hyperbolic) disk, such that type I vertices are mapped to the interior, type II vertices to the boundary, and the vertex $1$ is mapped to $0$. 

The form that is integrated over is a product of one-forms $\nu(u,v)$, one for each edge $(u,v)$ present in the graph $\Gamma$. Concretely, the form $\nu(u,v)=d\kappa(u,v)$ is the differential of the (hyperbolic) angle $\kappa(u,v)$ between the (hyperbolic) straight lines connecting $u$ and $1$ and $u$ and $v$. 


\subsection{Nomenclature and configuration space}
We do not use the standard numberings $1,\dots, m$ and $\bar{1},\dots,\bar{n}$ for the vertices, but instead take the ``customised'' names indicated in the picture (Fig. \ref{fig:theLiegraph}) for convenience. We will take the same names for the appropriate coordinate functions on the configuration space. The vertex $z$ plays the role of the vertex ``1'' and is fixed at 0. So e.g., $z=0$, $|w|,|a_j|,|b|<1$, and $U,V\in S^1$. Concretely, denote by $\alpha, \beta\in [0,1)$ the (normalized) arguments of $U,V$, i.e., $U=e^{2\pi i\alpha},V=e^{2\pi i\beta}$.

\subsection{The polylogarithm functions}
Denote the polylogarithm functions by $\li_n(x)$. They are defined on the unit disk $D:=\{|x|<1\}\subset \C$ through their Taylor series 
\[
\li_n(x) = \sum_{j=1}^\infty \frac{x^j}{j^n}
\]
and can be analytically continued to a multi-valued function on the complex plane. For us, only the continuation to the closed unit disk $\bar{D}=\{|x|\leq 1\}\subset \C$, and only integer $n=0,1,..$ are important. The function $\li_0(x)=\frac{x}{1-x}$ can be extended to $\bar{D}\setminus \{1\}$, having a pole at $x=1$. The function $\li_1(x)=-\log(1-x)$ can be extended to $\bar{D}\setminus \{1\}$, having a branch point at $x=1$. For $n>1$ the Taylor series converges on $\bar{D}$ and defines a continuous function. In fact, the limiting function $y\mapsto \li_n(e^{2\pi i y})$ is smooth on $(0,1)$.

The polylogarithm functions satisfy the identities
\begin{align*}
x \frac{d\li_n(x)}{d x} &= \li_{n-1}(x) \\
\li_{n}(e^{2\pi i y})+(-1)^n\li_{n}(e^{-2\pi i y}) &= -\frac{(2\pi i)^n}{n!}B_n(y)
\end{align*}
where $n=0,1,2,\dots$, $y\in (0,1)$ and $B_n(y)$ are the Bernoulli polynomials (see the appendix for a proof).

\section{The hyperbolic angle in terms of polylogarithms}
We will need the explicit form of the hyperbolic angle function $\kappa(w,U)$ corresponding to the edge $(w,U)$. It can be written down using polylogarithms as follows
\[
\kappa(w,U) = \alpha - \frac{1}{\pi} \Im\left( \li_1(w\bar{U}) - \li_1(w) \right).
\]
Similarly, 
\[
\kappa(w,V) = \beta - \frac{1}{\pi} \Im\left( \li_1(w\bar{V}) - \li_1(w) \right).
\]

\section{The strategy}
Our goal is to compute the Kontsevich integral (i.e., the weight) of the graph shown in Figure \ref{fig:theLiegraph}. Let us use the the following strategy:

\begin{enumerate}

\item Integrate out the variables $a_1,..,a_4$. This yields a function $F(U)$. It is a result due to C. Torossian  (see the appendix) that $F$ is a Bernoulli polynomial in $\alpha$, which we write as the real part of a polylogarithm function: $F(U)\propto \frac{1}{\pi^4}\li_4(U)+c.c.$. Here and in the following ``$c.c.$'' denotes the complex conjugate of the preceding expression. 

\item Next we integrate out the vertices $w$ and $b$, yielding a form $G(U,V)d\beta+(\text{irrelevant})d\alpha$. This is the major computation. The function $G(U,V)$ will be given as a Laurent series in $U,V$. Actually, it will be sufficient to compute the $V^0$-terms, which we denote by $G_0(U)$.

\item The remaining integration over $\alpha$ and $\beta$ can be simplified by the following observation: We claim that one can integrate over $(\alpha,\beta)\in (0,1)\times (0,1)$ and divide by 2 afterwards, instead of integrating over $0<\alpha<\beta<1$. The proof is by a reflection argument.\footnote{Concretely, the reflection wrt. the $x$-axis (i.e., complex conjugation) maps the region $\{\alpha > \beta\}$ to the original domain of the integral $\{\alpha < \beta\}$. The orientation is \emph{not} changed: Both $\alpha$ and $\beta$ contribute a $-1$ each, and each type I vertex except the fixed one contributes another $-1$. Hence the orientation changes by a factor $(-1)^{2+m-1}=1$. Pulling back the weight form, one obtains an additional $-1$ for each edge, hence an overall $+1$ since the number of edges is even.}

\item The integral we are left with is 
\[
\int_0^{1} d\alpha \int_0^{1} d\beta F(U)G(U,V)
\]
The integrals pick out the $U^0V^0$-term in the Laurent series of $F(U)G(U,V)$. Since $F(U)$ contains no $V$-terms, it is sufficient to compute the part $G_0(U)$ of $G(U,V)$, that consists of the $V^0$-terms. 

\end{enumerate}

\subsection{Laziness policy}
If one wants to compute the weight exactly, i.e., including all rational summands, one ends up managing expressions with roughly 50 terms. To avoid this, we compute everything only modulo uninteresting, but computable, rational prefactors and summands. In particular, note that $F(U)$ is a polynomial with rational coefficients in $\alpha$. Hence it is sufficient to determine $G(U,V)$ up to a rational-coefficient polynomial in $\alpha, \beta$, because the integral over a rational-coefficient polynomial yields a rational number.

We use the notation $A \sim B$, meaning that $A=p B+q$ for some $0\neq p\in \mathbb{Q}$, and some rational coefficient polynomial $q$ in $\alpha, \beta$. 

\section{The main computation}
Here the above plan is carried out, i.e., the integration over $b$ and $w$ is performed.

\subsection{Integration over $b$}
\label{sec:bint}
This integration has been done in the appendix, we quote the result:
\begin{lemma}
The $b$-integration yields a form (up to a rational prefactor) 
\[
\left(\kappa(w,V) - H(\beta-\alpha)
\right) d\kappa(w,U) d\kappa(w,V)
\]
where $H$ is the Heaviside function, 
\[
H(x)=\begin{cases}
0 &\quad \text{for $x\leq 0$}\\
1 &\quad \text{for $x>0$}
\end{cases}.
\]
\end{lemma}

The ``$-H(\dots)$'' is irrelevant, since, in effect, it contributes only a rational summand, due to the following Proposition.
\begin{prop}
The above graph, but with $b$ omitted and an extra edge $(w,U)$ added has rational weight.
\end{prop}
\begin{proof} 
See Corollary \ref{cor:liegrQ} in Appendix \ref{sec:cycinv}.
\end{proof} 

Hence the integral we need to compute is 
\[
\eta = \int_w \kappa(w,V) d\kappa(w,U)  d\kappa(w,V) d\phi.
\]

Here, the $\phi$ is the argument of $w$, i.e., $w=r e^{2\pi i\phi}$, and the form $d\phi$ is just the Kontsevich angle form corresponding to the edge $(z,w)$. (remember that $z=0$).

The form $\eta$ has a $d\alpha$-component and a $d\beta$-component, but for us only the former one is relevant (... since $z$ and $V$ are connected by an edge).

So we need to compute the integral
\[
\eta_\beta = \left( \int_0^1 dr \int_0^{1} d\phi
\kappa(w,V) \p_\alpha \kappa(w,U)  \p_r \kappa(w,V)
\right) d\alpha.
\]
The function in the bracket is just the $G(U,V)$ of the previous section.

\subsection{Computing the derivatives}
The polylogarithm satisfies $x \p_x \li_n{x}=\li_{n-1}(x)$. Hence we see that 
\[
\p_r \kappa(w,V) = -\frac{1}{\pi r} \Im\left( \li_0(w\bar{V}) - \li_0(w) \right)
\]
and
\[
\p_\alpha \kappa(w,U) = 1 + 2 \Re\left( \li_0(w\bar{U})  \right).
\] 
Putting everything together, the integral we need to compute becomes
\begin{multline}
\label{equ:GofUV}
G(U,V) \sim  \int_0^1 \frac{dr}{\pi r} \int_0^{1} d\phi
\left( \beta - \frac{1}{\pi} \Im\left( \li_1(re^{2\pi i\phi}\bar{V})  -  \li_1(re^{2\pi i\phi})  \right) \right)
\\ 
\Im\left( \li_0(re^{2\pi i\phi}\bar{V}) - \li_0(re^{2\pi i\phi}) \right)
\left(1 +2 \Re \li_0(re^{2\pi i\phi}\bar{U})  \right).
\end{multline}

We next want to simplify this expressions by identifying and throwing away irrelevant terms.

\subsection{First irrelevant term: $\beta$} 
The $\beta$ in the leftmost bracket is irrelevant, since it contributes a polynomial in $\alpha,\beta$ with rational coefficients to $G(U,V)$. This follows from the following Lemma.
\begin{lemma}
\label{lem:UV}
Let $U=e^{2\pi i\alpha}$ and $V=e^{2\pi i\beta}$, with $0\leq \alpha, \beta< 1$ and $\alpha\neq \beta$. Then 
\begin{align*}
&\frac{2}{(2\pi i)^{m+n+1}}\int_0^1 \frac{dr}{r} \int_0^{1} d\phi
(\li_m(re^{2\pi i\phi}\bar{V})-(-1)^m c.c.) ( \li_m(re^{2\pi i\phi}\bar{U}) +(-1)^n c.c.)
\\
&\quad \quad = 
-\frac{(-1)^n}{(n+m+1)!}
B_{m+n+1}(\alpha-\beta+H(\beta-\alpha) )
\sim 0
\end{align*}

\end{lemma} 
\begin{proof}
Taylor expansion of the integrand yields
\begin{align*}
&\int_0^1 \frac{dr}{r}
\int_0^{1} d\phi \sum_{j,k\geq 1} \frac{1}{j^m k^n} r^{j+k} e^{2\pi ij\phi}\bar{V}^j (e^{2\pi ik\phi}\bar{U}^k+(-1)^n e^{-2\pi ik\phi}U^k) - (-1)^{n+m} c.c.
\\ \quad\quad &= \lim_{\lambda\uparrow 1}
\int_0^\lambda \frac{dr}{r} \int_0^{1} d\phi (\dots)
\\ \quad\quad &=(-1)^n \lim_{\lambda\uparrow 1} \sum_j \frac{\lambda^{2j}}{2j^{m+n+1}} (\bar{V}U)^j -(-1)^{n+m}c.c
\\ \quad\quad &= \frac{(-1)^n}{2} \left( \li_{m+n+1}(\bar{V}U)+ (-1)^{m+n+1}c.c\right) 
\\ \quad\quad &=
 -\frac{(-1)^n}{2} \frac{(2\pi i)^{m+n+1}}{(n+m+1)!} B_{m+n+1}(\alpha-\beta + H(\beta-\alpha) )
\end{align*}
Here the regularization with the $\lambda$ is needed since the convergence radius of the polylogarithm's Taylor series is 1. The interchange of $r$-integral and sum would otherwise not be justified. 
\end{proof}

\subsection{Second irrelevant term: 1}
The summand ``1'' occuring in the rightmost brackets in \eref{equ:GofUV} will contribute only a polynomial in $\alpha$, $\beta$. The proof of this claim is almost the same as the one in the last section, using Lemma \ref{lem:UV} with $m=0,n=1$.

\subsection{Third irrelevant term: Everything containing $V$}
As said above, we are only interested in those terms of $G(U,V)$ that contain no $V$ or $\bar{V}$. When expanding the polylogarithms into Taylor series and performing the $\phi$-integral, a generic term will look like this:
\[
\sum_{j,k,l\geq 1} (\cdots)\cdot U^{\pm j} V^{\pm k} V^{\pm l} \delta(\pm j\pm k \pm l)
\] 
where the signs in front of the two $j$'s occuring are equal, as are those in front of the $k$'s and $l$'s, yielding 8 possible sign combinations. One can see that for neither of them the $V$'s cancel, since this would require $j=0$.

Hence we can kill all terms contributing nonzero powers of $V$ by simply setting $V=\bar{V}=0$ and using that $\li_n(0)=0$. We arrive at 

\[
G_0(U) \sim \frac{1}{\pi^2} \int_0^1 \frac{dr}{r} \int_0^1 d\phi
\Im\li_1(re^{2\pi i\phi})
\Im\li_0(re^{2\pi i\phi})
\Re\li_0(re^{2\pi i\phi}\bar{U}).
\]

\subsection{Computation of the integral}
In the above form, the integral is easily computed. However, as in the proof of Lemma \ref{lem:UV}, we need to introduce a cutoff $\lambda$ on the $r$-integral to justify the interchange of integrals and sums. For now, we keep the cutoff-dependence explicit, and take $\lambda\to 1$ at the end of the calculation. 

\begin{align*}
G_0(U, \lambda) &\sim \frac{1}{\pi^2} \sum_{j,k,l\geq 1} \int_0^\lambda \frac{dr}{r} r^{j+k+l} \int_0^{1} d\phi
  \frac{e^{2\pi ij\phi}}{j} 
  (e^{2\pi ik\phi}- e^{-2\pi ik\phi})
  (e^{2\pi il\phi}\bar{U}^l+e^{-2\pi il\phi}U^l) +c.c.
\\ &= \frac{1}{\pi^2}
\sum_{j,k,l\geq 1}  \frac{\lambda^{j+k+l}}{j(j+k+l)} (U^l \delta(j+k-l) - \bar{U}^l \delta(j-k+l) - U^l \delta(j-k-l)) +c.c.
\\ &= 
\frac{1}{\pi^2}\left(
\sum_{j,k\geq 1} \frac{\lambda^{2(j+k)}}{2j(j+k)} U^{j+k} - \sum_{j,l\geq 1}  \frac{\lambda^{2(j+l)}}{2j(j+l)}\bar{U}^l - \sum_{k,l\geq 1}  \frac{\lambda^{2(k+l)}}{2(k+l)^2} U^l 
\right) + c.c.
\\ &= 
\frac{1}{2\pi^2}\sum_{j,k\geq 1} \lambda^{2(j+k)} \left( \frac{U^{j+k}}{j(j+k)} - \frac{\bar{U}^{k}}{j(j+k)} - \frac{U^{j}}{(j+k)^2} \right) + c.c.
\end{align*}
In the last line, we relabeled the summation variables.

\section{Putting everything together}
We next evaluate
\[
\frac{1}{\pi^4} \int_0^{1}d\alpha \int_0^{1}d\beta
\Re\li_4(U)
G_0(U,\lambda)
\sim
\frac{1}{\pi^6} \sum_{j,k\geq 1} \lambda^{2(j+k)} \left( \frac{1}{j(j+k)^5} - \frac{1}{jk^4(j+k)} - \frac{1}{j^4(j+k)^2} \right)
\]
To evaluate these sums, introduce the notation 
\[
H_{n,r}= \sum_{j=1}^n \frac{1}{j^r}
\]
for the (generalized) harmonic numbers. Note that $\sum_{j\geq 1}\frac{1}{j(j+k)}=H_{k,1}/k$ (see Lemma \ref{lem:jjpn} in the Appendix). Then our integral becomes 
\[
\frac{1}{\pi^6}\left( \sum_{n\geq 1} \frac{H_{n-1,1}}{n^5} - \sum_{n\geq 1} \frac{H_{n,1}}{n^5} - \sum_{n\geq 1} \frac{H_{n-1,4}}{n^2} \right)
=
-\frac{1}{\pi^6}\zeta(6) - \frac{1}{\pi^6} \sum_{n\geq 1} \frac{H_{n-1,4}}{n^2}
\]
where we have removed the cutoff. The first term is rational. The second term is computed in Lemma \ref{lem:zetasums} in the appendix. The result is that the weight of the graph is
\[
\sim \frac{\zeta(3)^2}{\pi^6}.
\]

\appendix

\section{Polylogarithm and Bernoulli Polynomials}
\begin{lemma}
For $x\in (0,1)$ and $n=0,1,..$
\begin{align*}
B_n(x) = -\frac{n!}{(2\pi i)^n}\left( \li_{n}(e^{2\pi i x})+(-1)^n\li_{n}(e^{-2\pi i x}) \right).
\end{align*}

\end{lemma}
\begin{proof}
The Bernoulli polynomials are defined recursively by the equations $B_n'=nB_{n-1}$ and $\int_0^1B_n = \delta_{n0}$. The first one is easily checked, as is the second one for $n>1$. For $n=0$, explicitly
\[
-\left( \li_{0}(e^{2\pi i x})+\li_{0}(e^{-2\pi i x}) \right)
=
-\frac{e^{2\pi i x}}{1-e^{2\pi i x}}-\frac{1}{e^{2\pi i x}-1} = 1 = B_0.
\]
For $n=1$:
\[
-\frac{1}{2\pi i}\left( \li_{1}(e^{2\pi i x})-\li_{1}(e^{-2\pi i x}) \right)
=
\frac{1}{\pi}\arg (1-e^{2\pi i x})
=-\frac{1}{\pi} \frac{\pi-2\pi x}{2} = x-\frac{1}{2} = B_1(x)
\]
\end{proof}

\section{Bernoulli graphs}
\label{sec:berngraphs}
In this section it is shown that the following graph yields a Bernoulli polynomial in $x$. This computation is due to Charles Torossian.

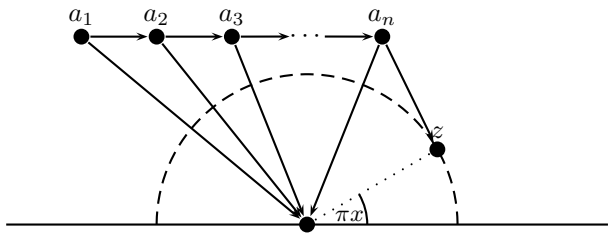
\begin{center}
\psset{unit=1cm}
\psset{arrowscale=1}

\begin{pspicture}(-4,-1)(4,3)
  \cnode*(-3,2.5){3pt}{a1}
  \cnode*(-2,2.5){3pt}{a2}
  \cnode*(-1,2.5){3pt}{a3}  
  \cnode*(1,2.5){3pt}{an}
  \rput[c](0,2.5){\rnode{dots}{$\cdots$}}
  \cnode*(2;30){3pt}{z}  

  \cnode*(0,0){3pt}{l}

  \psline{-}(-4,0)(4,0)
  \psarc[linestyle=dashed]{-}(0,0){2}{0}{180}
  \psarc{-}(0,0){.8}{0}{30}
  \ncline[linestyle=dotted]{l}{z}

  \psset{arrows=->}
  \ncline{a1}{a2}
  \ncline{a2}{a3}

  \ncline{an}{z}
  
  \ncline{a1}{l}
  \ncline{a2}{l}
  \ncline{a3}{l}
  \ncline{a3}{dots}
  \ncline{dots}{an}
  \ncline{an}{l}

  \uput{1ex}[90](a1){$a_1$}
  \uput{1ex}[90](a2){$a_2$}  
  \uput{1ex}[90](a3){$a_3$}
  \uput{1ex}[90](an){$a_n$}  
  \uput{1ex}[90](z){$z$}
  \uput{2.5ex}[15](l){\small $\pi x$}

\end{pspicture}

\end{center}

Let $\Gamma_n(x)$ be the value of the graph for fixed base angle $\pi x$. The integral $\int_0^1 \Gamma_n(x)dx$ vanishes for $n>0$ by a Kontsevich Lemma and is 1 for $n=0$. The differential of $\Gamma_n(x)$ is represented by the sum of all contractions of the graph. But the only nontrivial contraction possible is the contraction of the edge $(a_n,z)$. This yields the form $\Gamma_{n-1}(x) dx$. Hence $\Gamma_n' = \Gamma_{n-1}(x)$ and we obtain
\[
\Gamma_n(x) = \frac{1}{n!} B_n(x).
\]

\section{The $b$-integration}
Here we compute the integral over the position of vertex $b$ as promised in section \ref{sec:bint}. To restate the problem, we integrate over $b$ the weight form given by the following graph: 

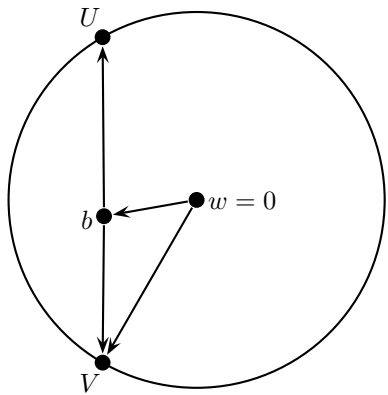
\begin{center}
\psset{arrowscale=1.5}
\begin{pspicture}(-3,-3)(3,3)
  \pscircle(0,0){2.5}
  \cnode*(0,0){3pt}{w}
  \cnode*(2.5;120){3pt}{U}
  \cnode*(2.5;-120){3pt}{V}
  \cnode*(1.25;190){3pt}{b}
  
  \psset{arrows=->}
  \ncline{w}{V}
  \ncline{b}{V}
  \ncline{b}{U}
  \ncline{w}{b}
  
  \uput{1ex}[120](U){$U$}
  \uput{1ex}[-120](V){$V$}  
  \uput{1ex}[0](w){$w=0$}  
  \uput{1ex}[190](b){$b$}
\end{pspicture}
\end{center}

We use the hyperbolic symmetries to put the vertex $w$ at $0$ for simplicity. 

Introduce polar coordinates $b=re^{2\pi i\phi}$ and $U=e^{2\pi i \alpha},V=e^{2\pi i \beta}$. 
Note that the result of the integration will be a two-form in $\alpha$ and $\beta$, i.e., be of the form $f(U,V)d\alpha d\beta$. The edge $(w,V)$ contributes a factor $d\beta$. The edge $(b,U)$ contributes a form $(\dots)d\alpha +(\dots)dr+(\dots)d\phi$. However, since it is the only edge that can contribute a $d\alpha$-term, the only relevent factor is
\[
(\dots)d\alpha = \p_\alpha ( \alpha - \frac{1}{\pi}\Im (\li_1(b \bar{U})-\li_1(b)) ) d\alpha =  ( 1 + 2\Re \li_0(b \bar{U}) ) d\alpha.
\]

The edge $(w,b)$ contributes a $d\phi$. Hence the edge $(b,V)$ has to contribute the missing factor of $dr$, namely
\[
\p_r ( \beta - \frac{1}{\pi}\Im (\li_1(b \bar{V})-\li_1(b)) ) dr =  -( \Im \li_0(b \bar{V}) - \Im \li_0(b)) \frac{dr}{\pi r}.
\]

The function $f(U,V)$ is hence given by:
\begin{align*}
f(U,V)
&=
-\int_0^1 \frac{dr}{\pi r} \int_0^1 d\phi  ( 1 +2 \Re \li_0(b \bar{U}) )( \Im \li_0(b \bar{V}) - \Im \li_0(b))
\\ &=
-2\int_0^1 \frac{dr}{\pi r} \int_0^1 d\phi  ( \Re \li_0(b \bar{U}) )( \Im \li_0(b \bar{V}) - \Im \li_0(b))
\end{align*}
For the last line we used that in the Taylor expansion of the right hand bracket, there is no constant term in $\phi$, hence the $\phi$-integral over it will vanish. Compute:
\begin{align*}
g(U,V) &:=
\int_0^1 \frac{dr}{r} \int_0^1 d\phi   \Re \li_0(b \bar{U}) \Im \li_0(b \bar{V}) 
\\ &=
- \frac{1}{4i} \left( \int_0^1 \frac{dr}{r} \sum_{j,k\geq 1} r^{j+k} \bar{U}^j V^k \delta(j-k) - c.c \right)
\\ &=
- \frac{1}{8i} \left( \sum_{j \geq 1} \frac{ (\bar{U}V)^j }{j} - c.c \right) = - \frac{1}{4} \Im \li_1(\bar{U}V)
\end{align*}

Hence we obtain
\begin{align*}
f(U,V) 
&=
- \frac{2}{\pi}(g(U,V)-g(U,1))
=
-\frac{1}{2\pi}\left(-\Im \li_1(\bar{U}V)+\Im \li_1(\bar{U}) \right)
\\ 
&= 
\begin{cases}
-\frac{1}{2}\left(B_1(\beta-\alpha)-B_1(1-\alpha) \right)
= \frac{1}{2}\left(-\beta + 1) \right) &\quad\quad \text{for $\beta>\alpha$} \\
-\frac{1}{2}\left(B_1(\beta-\alpha+1)-B_1(1-\alpha) \right)
= -\frac{1}{2}\beta &\quad\quad \text{for $\beta<\alpha$} 
\end{cases}
\end{align*}

\section{Some sums}
We collect here some helper results concerning (multiple) $\zeta$-function related sums.

\begin{lemma}
\label{lem:jjpn}
For $n$ a nonnegative integer, we have
\begin{align*}
\sum_{j\geq 1} \frac{n}{j(n+j)} &= H_{n,1} \\
\sum_{j\geq 1, j\neq n} \frac{n}{j(n-j)} &= H_{n,1}-\frac{2}{n}
\end{align*}
\end{lemma}
\begin{proof}
\begin{align*}
\sum_{j\geq 1} \frac{n}{j(n+j)} 
&= \lim_{\lambda\uparrow 1 } \sum_{j\geq 1} \frac{n\lambda^j}{j(n+j)}
= \lim_{\lambda\uparrow 1 }\left( \sum_{j\geq 1} \frac{\lambda^j}{j}- \sum_{j\geq 1} \frac{\lambda^j}{n+j}\right)
\\ &= \lim_{\lambda\uparrow 1 }\left( (1-\lambda^{-n}) \log(1-\lambda)  +  \sum_{j=1}^n \frac{\lambda^(j-n)}{j}\right) = 0 + H_{n,1}
\end{align*}

For the second equality, note that
\begin{align*}
\sum_{j\geq 1, j\neq n} \frac{n}{j(n-j)} 
=
-\sum_{k\geq 1} \frac{n}{(n+k)k}+ \sum_{j= 1}^{n-1}(\frac{1}{j} +\frac{1}{n-j})
=-H_{n,1} +2H_{n-1,1} = H_{n,1}-\frac{2}{n}
\end{align*}

\end{proof}

The following statements can essentially be found on the MathWorld-Homepage \cite{mworld}. 
\begin{lemma}
\label{lem:zetasums}
\begin{align*}
\sum_{n\geq 1} \frac{H_{n,4}}{n^2} &= \frac{25}{3}\zeta(6)-3\zeta(2)\zeta(4)-\zeta(3)^2\\
\sum_{n\geq 1} \frac{H_{n-1,4}}{n^2} &= \zeta(2,4)= \frac{22}{3}\zeta(6)-3\zeta(2)\zeta(4)-\zeta(3)^2 \\
2 \sum_{n\geq 1} \frac{H_{n,r}}{n^r} &= \zeta(r,1)+\zeta(r+1)=\zeta(r)^2-\zeta(2r) \\
2 \sum_{n\geq 1} \frac{H_{n,1}}{n^m} &= \zeta(m,1)+\zeta(m+1)=
(m+2) \zeta(m+1) - \sum_{n=1}^{m-2} \zeta(m-n)\zeta(n+1)
\end{align*}
\end{lemma}
\begin{proof}
The first statement is a direct consequence of the second. The last two statemets are needed to prove the second. Let us start with the third:
\[
2 \sum_n \frac{H_{n,r}}{n^r} = 2\sum_{m<n} \frac{1}{m^rn^r}
= \sum_{m\neq n} \frac{1}{m^rn^r} = \sum_{m, n} \frac{1}{m^rn^r}- \sum_{m= n} \frac{1}{m^rn^r}=\zeta(r)^2-\zeta(2r)
\]

Next consider the fourth:
\begin{align*}
\sum_{n=1}^{m-2} \zeta(m-n) \zeta(n+1) - (m-2) \zeta(m+1) 
&=
\sum_{i\neq j} \sum_{n=1}^{m-2} \frac{1}{i^{n+1}j^{m-n}}
\\ &= 
-\sum_{i\neq j} \frac{j^{2-m}-i^{2-m}}{ij(j-i)}
\\ &= 
-2\sum_{i\neq j} \frac{1}{ij^{m-1}(j-i)}
\\ &= 
-2\sum_{j} \frac{1}{j^{m-1}} \sum_{i\geq 1, i\neq j} \frac{1}{i(j-i)}
\\ &= 
-2\sum_{j} \frac{H_{n,1}}{j^{m-1}} + 4\zeta(m+1)
\end{align*}

For the first statement of the Lemma, we use the following identity, coming from the partial fractions expansion of $\frac{1}{j^3(j+k)^3}$:
\begin{align*}
2\sum_{j,k}\frac{1}{j^3(j+k)^3}
&= 
\sum_{j,k}\frac{1}{j^3 k^3}-\frac{3}{k^4 j^2}+\frac{6}{k^5 j}-\frac{3}{k^4 (j+k)^2}-\frac{6}{k^5 (j+k)}
\\ &=
\zeta(3)^2- 3\zeta(4)\zeta(2) +6 \sum_{j,k}\frac{1}{k^4 j (j+k)} - 3 \sum_n\frac{H_{n-1,4}}{n^2}
\\ &=
\zeta(3)^2- 3\zeta(4)\zeta(2) +6 \sum_{k}\frac{H_{k,1}}{k^5} - 3 \sum_n\frac{H_{n-1,4}}{n^2}
\end{align*}
Using the third and fourth identities from the lemma, we obtain 
\begin{align*}
3 \sum_n\frac{H_{n-1,4}}{n^2}
&=
\zeta(3)^2 -3\zeta(4)\zeta(2)
+ 3\left(
7\zeta(6) - 2\zeta(4)\zeta(2) -\zeta(3)^2
\right)
- \zeta(3)^2+\zeta(6)
\\ &=
22\zeta(6)-9\zeta(4)\zeta(2)-3\zeta(3)^2
\end{align*}
\end{proof}

\section{Cyclic Invariance of Kontsevich's morphism}
\label{sec:cycinv}
At two points in this paper, relations between weights of graphs due to the cyclic invariance of Kontsevich's morphism are used. The underlying Theorem is the following.

\begin{prop}
Let $\Gamma$ be a Kontsevich graph. Let $\sigma \Gamma$ be the graph obtained by cyclically rotating the (labels on the) type II vertices $\bar{1},\dots,\bar{n}$. Then the weight $w_\Gamma$ satisfies
\[
(-1)^n w_{\sigma \Gamma} = \sum_{\Gamma'} (-1)^{n_1(\Gamma')}{w_{\Gamma'}}.
\]
Here the sum on the right is over all $\Gamma'$ that can be obtained from $\Gamma$ by replacing an arbitrary number of edges $(v_1,v_2)$ of $\Gamma$ by edges $(v_1,\bar{1})$. The number $n_1(\Gamma')$ is the number of edges in $\Gamma'$ connecting to $\bar{1}$.
\end{prop}
This proposition is just the reformulation in terms of graph weights of the result of Felder and Shoikhet \cite{fsh} that the image of divergence free polyvector fields under Kontsevich's morphism is cyclic.

We only need the following two Corollaries:
\begin{cor}
\label{cor:liegrQ}
Any Lie graph with exactly two edges hitting a type II vertex has weight $\frac{B_p B_q}{2}$ for some integers $p,q$.
\end{cor}
\begin{cor}
\label{cor:wheelweight}
The middle graph of Figure \ref{fig:stargraphs} has weight $\frac{B_k}{2\cdot k!}$.
\end{cor}

\section{Kontsevich and CBH product on duals of Lie algebras}
\label{sec:kontscbh}
The goal of this section is to give a short proof of the following Theorem, mentioned in the introduction.
\begin{thm}
The CBH product is isomorphic to the Kontsevich product on duals of Lie algebras. The automorphism of $S\alg{g}$ mapping the two star products onto each other is given by the Duflo map
\[
\exp\left( \sum_{j\geq 1}\frac{B_{2j}}{4j (2j)!} \tr(\ad_\p^{2j}) \right).
\]
\end{thm}
We want to emphasize that the result is well known. It follows from results of Kontsevich (see \cite{kontsevich}, section 8) and Shoikhet \cite{shwheel}. For nilpotent Lie algebras, it has been explicitly proven by V. Kathotia \cite{kathotia}.

Denote by $\star_{CBH}$ the CBH product and by $\star_{DK}$ the pullback of Kontsevich's product under the Duflo map. 
\begin{lemma}[First Reduction]
Let $\star$ be some associative product on $S\alg{g}$. If for any $X,Y\in \alg{g}$ and $n=1,2,..$ we have that 
\[
X^n \star_{CBH} Y = X^n \star Y
\]
then $\star = \star_{CBH}$.
\end{lemma}
\begin{proof}
We have to show that $P\star_{CBH} Q = P \star Q$ for any polynomials $P,Q$ on $\alg{g}^*$. By polarization, it is sufficient to show that $X^n \star_{CBH} Y^m = X^n \star Y^m$ for any $X,Y$ as above and $n,m=1,2,\dots$ Also by polarization, from the assumption in the lemma it follows that $P \star_{CBH} Y = P \star Y$ for any polynomial $P$. Using associativity we then compute
\begin{multline*}
X^n \star_{CBH} Y^m = X^n \star_{CBH} Y^{\star_{CBH} m}
=(\cdots(X^n \star_{CBH} Y) \star_{CBH} \cdots )\star_{CBH} Y
=\\=(\cdots(X^n \star Y) \star \cdots \star Y)
= X^n \star Y^{\star m}.
\end{multline*}
The result follows by the observation 
\[
Y^{\star m} = (\cdots(Y\star Y)\star \cdots )\star Y= Y^{\star_{CBH} m}= Y^m.
\]
\end{proof}

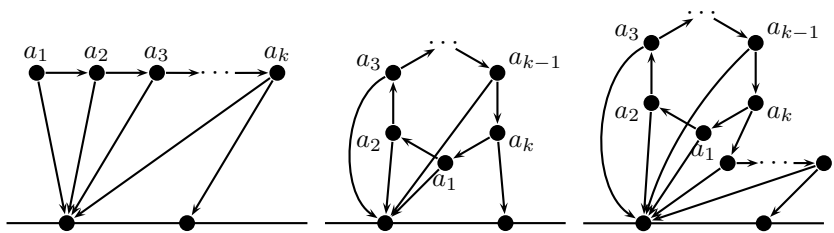
\begin{figure}
\psset{unit=.8cm}
\psset{arrowscale=1}
\begin{pspicture}(-2,-.5)(3,3)
  \cnode*(-1.5,2.5){3pt}{a1}
  \cnode*(-0.5,2.5){3pt}{a2}
  \cnode*(0.5,2.5){3pt}{a3}  
  \cnode*(2.5,2.5){3pt}{an}
  \rput[c](1.5,2.5){\rnode{dots}{$\cdots$}}
  
  \cnode*(1,0){3pt}{r}  
  \cnode*(-1,0){3pt}{l}

  \psline{-}(-2,0)(3,0)

  \psset{arrows=->}
  \ncline{a1}{a2}
  \ncline{a2}{a3}

  \ncline{an}{r}
  
  \ncline{a1}{l}
  \ncline{a2}{l}
  \ncline{a3}{l}
  \ncline{a3}{dots}
  \ncline{dots}{an}
  \ncline{an}{l}

  \uput{1ex}[90](a1){$a_1$}
  \uput{1ex}[90](a2){$a_2$}  
  \uput{1ex}[90](a3){$a_3$}
  \uput{1ex}[90](an){$a_k$}  
\end{pspicture}
\begin{pspicture}(-2,-2.5)(2,1.5)
  \cnode*(1;30){3pt}{ak1}
  \cnode*(1;-30){3pt}{ak}
  \cnode*(1;-90){3pt}{a1}
  \cnode*(1;-150){3pt}{a2}  
  \cnode*(1;-210){3pt}{a3}  
  \rput[c](1;90){\rnode{dots}{$\cdots$}}
  
  \cnode*(1,-2){3pt}{r}  
  \cnode*(-1,-2){3pt}{l}

  \psline{-}(-2,-2)(2,-2)

  \psset{arrows=->}
  \ncline{a1}{a2}
  \ncline{a2}{a3}
  \ncline{ak}{a1}
  \ncline{ak1}{ak}

  \ncline{ak}{r}
  
  \ncline{a1}{l}
  \ncline{a2}{l}
  \ncarc[arcangle=-60]{a3}{l}
  \ncline{ak1}{l}
  
  \ncline{a3}{dots}
  \ncline{dots}{ak1}

  \uput{1ex}[-90](a1){$a_1$}
  \uput{1ex}[-150](a2){$a_2$}  
  \uput{1ex}[-210](a3){$a_3$}
  \uput{1ex}[-30](ak){$a_k$}  
  \uput{1ex}[30](ak1){$a_{k-1}$} 
\end{pspicture}
\begin{pspicture}(-2,-3)(2,1.5)
  \cnode*(1;30){3pt}{ak1}
  \cnode*(1;-30){3pt}{ak}
  \cnode*(1;-90){3pt}{a1}
  \cnode*(1;-150){3pt}{a2}  
  \cnode*(1;-210){3pt}{a3}  
  \rput[c](1;90){\rnode{dots}{$\cdots$}}
  
  \cnode*(0.4,-1.5){3pt}{b1}  
  \rput[c](1.2,-1.5){\rnode{bdots}{$\cdots$}}
  \cnode*(2.0,-1.5){3pt}{br}  
  
  \cnode*(1,-2.5){3pt}{r}  
  \cnode*(-1,-2.5){3pt}{l}

  \psline{-}(-2,-2.5)(2,-2.5)

  \psset{arrows=->}
  \ncline{a1}{a2}
  \ncline{a2}{a3}
  \ncline{ak}{a1}
  \ncline{ak1}{ak}

  \ncline{ak}{b1}
  
  \ncline{a1}{l}
  \ncline{a2}{l}
  \ncarc[arcangle=-60]{a3}{l}
  \ncarc[arcangle=-15]{ak1}{l}
  \ncline{b1}{l}
  \ncline{br}{r}
  \ncline{br}{l}
  
  \ncline{a3}{dots}
  \ncline{dots}{ak1}
  \ncline{b1}{bdots}
  \ncline{bdots}{br}
  
  \uput{1ex}[-90](a1){$a_1$}
  \uput{1ex}[-150](a2){$a_2$}  
  \uput{1ex}[-210](a3){$a_3$}
  \uput{1ex}[-30](ak){$a_k$}  
  \uput{1ex}[30](ak1){$a_{k-1}$} 
\end{pspicture}
\caption{\label{fig:stargraphs} The three kinds of graphs occuring in the pullback of the Kontsevich star product by the Duflo map.}
\end{figure}

Now we can prove the Theorem. It is well known that $X^n \star_{CBH} Y = \sum_j {n \choose j} B_j X^{n-j} \ad_X^jY$. On the other hand, 
the graphs occuring in $X^n \star_{DK} Y$ are of one of the three types shown in Figure \ref{fig:stargraphs}. The left type has weight $\frac{B_k}{k!}$ by Appendix \ref{sec:berngraphs}, and yields the expression $X^n \star_{CBH} Y$. The right hand one yields no contribution since $\tr(\ad_X^r\ad_{\co{X}{Z}})=\tr(ad_X^r\co{\ad_X}{\ad_{Z}})=0$ for any $Z\in \alg{g}$ by cyclicity of the trace. The middle type comes in twice: Firstly from the twisting by the Duflo map with a weight $-\frac{B_k}{2 k \cdot k!}\cdot k=-\frac{B_k}{2\cdot k!}$. Secondly, from the Kontsevich formula, with a weight $\frac{B_k}{2\cdot k!}$ by Corollary \ref{cor:wheelweight}. So the graphs of the middle kind cancel and the Theorem has been shown.

\bibliographystyle{plain}
\bibliography{graphzeta}

\end{document}